\documentclass[11pt]{article}    % Specifies the document style.
\usepackage{fullpage,amsmath,graphicx,theorem}
\usepackage{subfigure}
\usepackage{amsfonts}
\usepackage{amssymb}
\usepackage{epsfig}
\newcommand {\beq}{\begin{equation}}
\newcommand {\eeq}{\end{equation}}
\newcommand {\bearn}{\begin{eqnarray*}}
\newcommand {\eearn}{\end{eqnarray*}}

\newtheorem{theorem}{Theorem}
\newtheorem{lemma}{Lemma}
\newtheorem{remark}{Remark}
\newtheorem{proposition}{Proposition}

\newtheorem{assumption}{Assumption}
\newtheorem{definition}{Definition}
\newtheorem{example}{Example}

\newenvironment{hangref}{\begin{list}{}{\setlength{\itemsep}{0pt}
\setlength{\parsep}{0pt}\setlength{\leftmargin}{+\parindent}
\setlength{\itemindent}{-\parindent}}}{\end{list}}

\title{Selecting the best system and multi-armed
                                                                          bandits}  % Declares the document's title.
\author{Peter Glynn \\
Stanford University \\
\\
Sandeep Juneja \\
 Tata Institute of Fundamental Research}

\begin{document}           % End of preamble and beginning of text.

\setlength {\baselineskip}{18 pt}

\maketitle                 % Produces the title.
\thispagestyle{empty}
\begin{abstract}
Consider the problem of finding a population or a probability distribution amongst many with the largest mean when these means are unknown but population samples can be simulated
or otherwise generated. Typically, by selecting largest sample mean population, it can be shown that false selection probability decays at an exponential rate. Lately, researchers have sought algorithms that guarantee that this probability is restricted to a small $\delta$  in order $\log(1/\delta)$  computational time by estimating the associated large deviations rate function via simulation. We show that such guarantees are misleading when populations have unbounded support even when these may be light-tailed. Specifically, we show that   any policy that identifies the correct population with probability at least $1-\delta$ for each problem instance requires 
infinite number of samples in expectation  in making such a determination in any problem instance. This suggests that some restrictions are essential on populations to devise $O(\log(1/\delta))$ algorithms with $1 - \delta$ correctness guarantees. We note that under restriction on population moments, such methods are easily designed, and that sequential methods from stochastic multi-armed bandit literature can be adapted to devise such algorithms.
\end{abstract}


%\pagebreak
\setcounter{page}{1}
\section{Introduction}
Suppose that we can sample independently from $d$ different probability distributions or populations $(X(i): i \leq d)$. Further, from each population $i$, we can generate independent identically distributed (iid) samples $(X(i,j): j \geq 1)$.  The distribution of $(X(i): i \leq d)$   is not known and our aim is to find the `best' population
\[
 i^* = \arg \max_{1 \leq j \leq d} EX(j).
\]
This and related problems of ascertaining the orders of expectations of populations are referred to as ordinal optimization in simulation literature.
The key observation is that ordinals are learnt correctly faster than the values of expectations. More specifically, when the
random variables $(X(i): i \leq d)$ are light-tailed (a distribution is said to be light-tailed if its   moment generating function is finite in a neighbourhood of zero) when $n$ samples are generated from each population and the population with the largest sample mean is selected
as the best one, the probability of false selection, $P(FS)$,
\[
P(\bar{X}_n(i^*) < \max_{j \ne i^*}\bar{X}_n(j))
\]
decays exponentially in $n$, where $\bar{X}_n(i) = \frac{1}{n} \sum_{k=1}^n X(i,k)$, while
it is well known through central limit theorem that the rate of convergence of
$\bar{X}_n(i)$  to $EX(i)$ is a much slower $n^{-1/2}$.

This exponential decay rate of false selection probability suggests that for small but positive $\delta$, one may be able to construct algorithms that, when acting on light-tailed distributions, by generating $O(\log(1/\delta))$  samples of $X(i)'s$, can make a correct selection with  the probability of false selection restricted to $\delta$.  In this paper we critically examine this proposition. We answer this in negative
 for algorithms acting on a collection of random variables with unbounded support on the right, under a mild technical condition. This includes algorithms 
 that act on the
 collection of all light-tailed distributions with unbounded support on the right. In fact, we conclude that  expected number of samples
 generated by such an algorithm would be infinite.   We provide a positive answer when further information on moments of the underlying random variables is available.

{\em Applications} of ordinal optimisation  in simulation arise in many set-ups including
in selecting a best design from a set of competing designs via simulation where  all of the designs may be modelled as discrete event dynamic systems.  Such systems include queueing systems, computer and communications
networks, manufacturing systems and transportation networks (see,
e.g., Ho, Srinivas and Vakili 1992 for some applications).

\subsection{Brief literature review}

\subsubsection{Simulation}

  Ho et. al. (1992) observed that determining ordinals amongst population means is faster than estimating the means.  Dai (1996) used large deviations theory to show  in a fairly general framework that for light-tailed  random variables the probability of false selection decays exponentially.
Chen, Lin, Yucesan and Chick (2000) considered the problem of ordinal optimization under the assumption that a fixed but large computation budget $n$ is available and the underlying random variables have a Gaussian distribution. They optimized the budget allocated to each population asymptotically as $n \rightarrow \infty$ so that the probability of false selection is minimized.  Glynn and Juneja (2004)  observed that
 if fraction $p_i n$ of computational budget is allocated to population $i$ for each $i$ ($p_i>0$, and $p_i n$ an integer), then the probability of false selection,
\begin{equation} \label{eqn:exp_bound101}
%P( \bar{X}_{i^*}(p_{i^*} n)  < \max_{i \ne i^*} \bar{X}_i(p_i n) ) \leq e^{-n H(p_1, \ldots, p_d)},
\lim_{n\rightarrow \infty} \frac{1}{n}\log P( \bar{X}_{i^*}(p_{i^*} n)  < \max_{i \ne i^*} \bar{X}_i(p_i n) )
= - H(p_1, \ldots, p_d),
\end{equation}
where $H(\cdot)$ is  concave.
 They, then maximized  $H(p_1, \ldots, p_d)$ under the constraint
$\sum_{i=1}^d p_i =1$ to determine the optimal allocations as $n \rightarrow \infty$ even for non-Gausssian distributions.  Significant literature since then has appeared that relies on large deviations analysis {(e.g., Hunter and Pasupathy 2013,  Szechtman and Yucesan 2008, Broadie, Han, Zeevi 2007, Blanchet, Liu, Zwart 2008, Shin, Broadie and Zeevi 2016).}

Substantial literature exists on selecting the best system amongst many alternatives using ranking/selection procedures, see
 Chernoff (1959), Paulson (1964), Bechhofer, Kiefer and Sobel (1968), Jennison, Johnstone and Turnbull (1982)  
 for some of the earlier references from the statistics literature.   
 
 See, e.g., Kim and Nelson (2001, 2003),
Branke, Chick and Schmidt (2007) for an overview, Frazier (2014) for a Bayesian inspired approach. Gaussian assumption is critical to most of the analysis
here. Typically,  an  `indifference-zone
formulation' is considered where it is assumed that there exists a known
$\epsilon>0$ such that
\begin{equation} \label{eqn:indiffzone}
 E X(i^*) \geq EX(j) +\epsilon
 \end{equation}
 for $j \ne i^*$ (see,
e.g., Nelson and Matejcik 1995).
 Kim and Nelson (2001) propose a number of algorithms for selecting the best system. They conduct asymptotic analysis
where the probability of false selection is kept fixed
at a specified $\delta>0$ while the indifference zone is set to $\epsilon$. They show that the proposed algorithms restrict $P(FS)$ to within $\delta$,
asymptotically as $\epsilon \rightarrow 0$, by generating  $O(\epsilon^{-2})$ samples.
In contrast, in this paper,  we primarily consider settings where the indifference zone is kept fixed and $\delta$ is either fixed, or it decreases to zero.

\subsubsection{Stochastic multi-armed bandits}

The problem of selecting the best system amongst many via sampling is  related to the vast, elegant and evolving literature  in learning and statistics community referred to as stochastic multi-armed bandit methods. Typically, this strand of  literature
refers to sampling from a population as `pulling an arm' and assumes that each such pull leads to a random reward whose distribution depends upon the population. The aim then is to develop optimal or near optimal sequential sampling strategy that
maximises the long term total expected reward or equivalently minimises the total expected regret (regret over $n$ trials is referred to as the reward that would have been realised in these trials if the `best' arm was pulled each time versus the actual reward realisation from a given sequential strategy - here best refers to the arm with the highest expected reward.
See, e.g.,  Bubeck and Cesa-Bianchi 2012, Cappe et. al. 2013; see Lai and Robbins 1985 for a seminal paper in this area). Optimal sequential strategies have to carefully manage the exploration-exploitation trade-offs in arm selections. Recently, these methods have also been used in the  {\em pure exploration setting} where the goal is to identify the arm with the highest expected reward in minimum
expected number of trials. See, e.g., Even-Dar, Mannor and Mansour (2002, 2006), Audibert and Bubeck (2010), Jamieson, Malloy, Nowak, \& Bubeck (2014),
Garivier and Kaufmann (2016) and Kaufmann, Cappe and Garivier (2016). 
 Thus, this problem is identical to the ordinal optimisation problem that we consider. A standard assumption both in minimising regret and the pure exploration settings is that the rewards from each arm are either Bernoulli, are bounded with known bounds, or have sub-Gaussian tails. 
 Garivier and Kaufmann (2016) and Kaufmann, Cappe and Garivier (2016) consider a somewhat general setting where the rewards come from a single-parameter
 exponential family.  In this setting, Garivier and Kaufmann (2016) develop lower bounds on expected number of samples from each arm,
 as well as algorithms that match these lower bounds asymptotically as $\delta \rightarrow 0$. 
  Bubeck, Cesa-Bianchi and Lugosi (2013), again consider the problem of minimising expected regret under the assumption that the rewards are unbounded but explicit bounds on their moments are known.

\subsection{Observations}

Suppose that, as in (\ref{eqn:exp_bound101}),  order $n$ samples are allocated to each population, and as is typically the case when the underlying random variables are
light-tailed, the probability of false selection,
$P(FS) \leq e^{-n I}$,  for some $I >0$.
If such an $I$ is known beforehand, setting
\[
n = - \frac{1}{{I}}\log(\delta) \mbox{ ensures  }  P(FS) \leq \delta.
 \]

Typically, such an $I$ is not known but needs to be estimated via simulation. One then hopes for algorithms that for $n= O(\log(1/\delta))$ ensure that  $P(FS) \leq \delta$, for each $\delta$, or this holds asymptotically, that is, 
\begin{equation} \label{eqn:efficiency}
\limsup_{\delta \rightarrow 0} P(FS) \delta^{-1} \leq 1,
\end{equation}
even when an indifference zone formulation is considered so that
there exists a  fixed and known $\epsilon>0$ and (\ref{eqn:indiffzone}) holds.

Note that
$O(\log(1/\delta))$ effort is necessary to achieve (\ref{eqn:efficiency}) in the sense that  if $\log(1/\delta)^{1-\tilde{\epsilon}}$ samples are generated, for any small $\tilde{\epsilon}>0$,
then for any $i$ and $A$ such that $P(X(i) \in A)>0$
\[
P(X(i) \in A)^{\log(1/\delta)^{1-\tilde{\epsilon}}}
=\delta^{\frac{\mbox{positive no.}}{\log(1/\delta)^{\tilde{\epsilon}}}} > \delta
\]
as $\delta \rightarrow 0$, so one cannot hope for (\ref{eqn:efficiency}) to hold.

Also observe that order $\log(1/\delta)^{1+\tilde{\epsilon}}$ is sufficient as for $c>0$
\[
\delta^{-1}P(FS) \leq  \delta^{-1} e^{-n I}= \delta^{-1}e^{-c  \log(1/\delta)^{1+\tilde{\epsilon}}I}
= \delta^{c \log(1/\delta)^{\tilde{\epsilon}}I-1}
\]
which goes to zero as $\delta \rightarrow 0$.

In this paper we shed light on conditions under which 
algorithms that require   $O(\log(1/\delta))$ samples and  satisfy  $P(FS) \leq \delta$ for each $\delta$, or more generally satisfy the relation
 (\ref{eqn:efficiency}), exist, as well as when such algorithms do not exist.

\subsection{Our contributions}

{\bf Our key result is negative} and is as follows -  Consider a collection of probability distributions
${\cal L}$ and  an algorithm that given any finite set of distributions belonging to  ${\cal L}$,
and a small $\delta>0$,
adaptively generates samples from these distributions and identifies the
distribution with the largest mean with the probability of false selection $P(FS) \leq \delta$.  These finite set of distributions are assumed to be well separated
in the sense that  the mean of the distribution with the highest mean exceeds the mean of other distributions by
 at least $\epsilon >0$. We also consider algorithms that only provide  asymptotic guarantees
 in the sense that
\begin{equation} \label{eqn:need1}
\limsup_{\delta \rightarrow 0} P(FS) \delta^{-1} \leq 1.
\end{equation}
We prove that if ${\cal L}$ comprises populations with unbounded support that are {\em KL right dense} (defined in Section 4), the expected number of samples so that $P(FS) \leq \delta$ for any  $\delta \in (0, 1/2)$, must be infinite. If the algorithm only provides asymptotic guarantees as in
(\ref{eqn:need1}),  then the expected number of samples cannot be
$O(\log(1/\delta))$.  Further, similar results
hold when the criteria for selecting the best design may not be the population mean but another function of the population distribution such as its specific quantile.

{\bf Our positive contributions} -
Under  explicitly available  upper bounds on convex, increasing functions
of underlying random variables, we develop random variable truncation as well as capping based
$O(\log(1/\delta))$ computation time algorithms that guarantee that $P(FS) \leq \delta$.
We also observe that  the sequential algorithms proposed in multi-armed bandit regret setting (see  Bubeck, Cesa-Bianchi and Lugosi  2013), when heavy tails are involved and  bounds on the $\alpha \in (1,2]$ moments of the underlying distributions are assumed known, are easily adapted to this {\em pure exploration setting} to provide
$O(\log(1/\delta))$ computation time sequential algorithms.  We also develop upper bounds on computational effort under these algorithms and suggest tweaks that lead to minor performance improvement guarantees.

Our attempt in presenting positive results is
to illustrate reasonable conditions under which $O(\log(1/\delta))$ computation time algorithms are feasible in principle. We do not attempt to excessively optimize the proposed implementations. A theory  involving general distributions that characterizes the  lower bounds on expected number of samples from each arm,
and develops algorithms that match these lower bounds, is an interesting area for future research. 
Recall that Garivier and Kaufmann (2016)  achieve this for single parameter exponential family of distributions.

\subsection{Roadmap}
 Section 2 contains our key negative result illustrating that, under mild regularity conditions, the algorithms that control the probability of false selection to within a small $\delta$ must require infinite number of samples in expectation.  In Sections 3 and 4 we provide some positive results. In Section 3,  we develop  $O(\log(1/\delta))$ algorithms when upper bounds on suitable moments of underlying random variables are available. In Section 4, we adapt recent results from multi-armed bandit related research  to our ordinal optimisation setting.  We end with a brief conclusion in Section 5.

All but the simpler  proofs are given in the appendix.

\section{The negative result}

%Let ${\cal L}$ denote a collection of probability distributions each with finite mean and with additional requirement  that given any $G \in {\cal L}$, there exists  $\tilde{G} \in {\cal L}$, whose mean is arbitrarily large and $H(G | \tilde{G}) < \infty$.

We consider an algorithm or a policy acting on any finite set of probability distributions taken from a collection 
${\cal L}$.
The policy does not know the underlying set of distributions but it can adaptively sample from them. The objective is to 
 to identify the distribution with largest mean while adaptively generating fewest number of samples, and providing guarantees
 on the  probability of error.

Specifically,  consider  a policy that acts on a set  of $d$ distributions $(F_i: i \leq d)$ in ${\cal L}$ . Let
  \begin{itemize}
  \item
$D_t$  denote the index of the distribution sampled from at stage $t$.
\item
$N_i(t)= \sum_{k=1}^t I(D_k=i)$ denote the number of samples
generated  from distribution $F_i$ by stage $t$.  
\item
  $X_{i,j}$ denote the sample $j$ generated from distribution
$F_i$. This is assumed to be independent of other previously generated samples. 
\item
${\cal F}_t$ denote the sigma-algebra associated with $(D_k, k \leq t)$ and $(X_{i,j}, i \leq d, j \leq N_i(t))$. 
\end{itemize}

Then, at stage $t+1$, the policy selects a distribution  $D_{t+1}$ to sample from 
where $D_{t+1}$ is measurable w.r.t.
${\cal F}_t$. Thus, if $D_{t+1}=s$, then the policy   generates a sample   $X_{s, N_{s}(t)+1}$ from distribution $F_{s}$.

Each policy considered is  further assumed to have  an almost sure finite  stopping time $\tau$ with respect  to the filteration $\{{\cal F}_t\}$.
At time $\tau$ the policy  announces
its decision  $J=i$ if it deems distribution $F_i$ to have the   largest mean.

\begin{definition} {
An algorithm or a  policy is referred to as  ${\cal P}(\epsilon, \delta)$, if  for any $\epsilon >0$ and $0 < \delta <1$,
given any finite collection of distributions from ${\cal L}$ such that the difference between the highest and the second highest mean amongst these distributions exceeds $\epsilon$, the policy, based on generated samples from these distributions,  identifies the correct distribution (the one with the largest mean) with the probability of false selection  $P(FS)$ guaranteed to be  $\leq \delta$.}
\end{definition}

\begin{definition} {
${\cal P}_{a}(\epsilon)$ denotes a larger class of policies that  differ from ${\cal P}(\epsilon, \delta)$ in that for any $\epsilon >0$, they only
guarantee that $P(FS)$ is asymptotically bounded from above by $\delta$.  That is,
\begin{equation} \label{eqn:asym_guarantee}
\limsup_{\delta \rightarrow 0} P(FS) \delta^{-1} \leq 1.
\end{equation}}
\end{definition}

Thus, roughly speaking, a policy ${\cal P}_{a}(\epsilon)$ is ${\cal P}(\epsilon, \delta)$ for all positive and sufficiently small $\delta$.
In many applications, it may be difficult to design a ${\cal P}(\epsilon, \delta)$ policy, however, designing
a ${\cal P}_{a}(\epsilon)$ policy may be easier, and often a practitioner may be satisfied with
asymptotic guarantees. (This is analogous to accepting central limit theorem based confidence intervals for population means
where the coverage guarantees are valid only asymptotically, see, e.g., Glynn and Whitt 1992). Note that if a policy is ${\cal P}(\epsilon, \delta)$
on ${\cal L}$, then it is also  ${\cal P}_{a}(\epsilon)$ on ${\cal L}$.

Suppose that the data was generated from  distributions
$(F_i: i \leq d)$ under probability measure $P$, and from
$(\tilde{F}_i: i \leq d)$ under probability measure $\tilde{P}$,
where  $F_i$ and $\tilde{F}_i$ are mutually absolutely continuous for each $i$.
Let $H(G | \tilde{G})= \int_{x \in \Re} \left ( \log \frac{d G}{d \tilde{G}}(x) \right ) d G(x)$ denote the
Kullback-Leibler distance between distributions $G$ and $\tilde{G}$.

Kaufmann, Cappe and Garivier (2016) Lemma 1, 
arrives at the following  elegant inequality (under the assumption that $\tau$ is a.s. finite under $P$),
\begin{equation} \label{eqn:Kauff1}
\sum_{i=1}^d E_P N_i(\tau) H(F_i | \tilde{F}_i)  \geq \sup_{{\cal E} \in {\cal F}_{\tau}} d(P({\cal E})|\tilde{P}({\cal E}))
\end{equation}
where for $p,q \in (0,1)$
\[
d(p, q) = p \log \frac{p}{q} + (1-p) \log \frac{1-p}{1-q},
\]
and ${\cal F}_{\tau}$ denotes the sigma algebra associated with stopping time $\tau$. 
This allows development of lower bounds $E_P N_i(\tau)$ under ${\cal P}(\epsilon, \delta)$
and  ${\cal P}_{a}(\epsilon)$ policies. Essential ideas for deriving such lower bounds
go back to  Lai and Robbins (1985) and Mannor and Tsitsiklis (2002). See Glynn and Juneja (2015), an earlier version of this paper,  for 
a derivation 
based on these ideas
when $\delta \rightarrow 0$.
Inequality  (\ref{eqn:Kauff1}) considerably simplifies
the development of lower bounds, and provides non-asymptotic bounds for any $\delta \in (0,1)$.

%{ Our key result of this section is that for any two distributions in ${\cal L}$
%with mean more than $\epsilon$ apart,  ${\cal P}_{asymp}(\epsilon, \delta)$ policy
%on  ${\cal L}$ requires in expectation  more than $O(\log(1/\delta))$ samp

\subsection{Lower bounds}  \label{secs:anal}

Our aim in this section is to show that ${\cal P}(\epsilon, \delta)$  policies
require infinite number of samples in expectation   to provide
the 
desired probabilistic guarantees, while ${\cal P}_{a}(\epsilon)$ policies
cannot provide the probabilistic guarantees  (\ref{eqn:asym_guarantee}) by generating only
$O(\log(1/\delta))$ samples when ${\cal L}$ comprises distributions with finite mean, unbounded support and additional mild restrictions.
 To illustrate this simply, we restrict our analysis
to $d=2$. Analysis extends easily  to $d >2$.

Consider distributions  $F$, $G$ and $\tilde{G}$ in ${\cal L}$
with respective means 
$\mu_F$,  $\mu_G$ and $\mu_{\tilde{G}}$
 For $\epsilon >0$,
  $\mu_G < \mu_F - \epsilon$
and $\mu_{\tilde{G}}> \mu_F + \epsilon$.
Further, 
the Kullback-Leibler distance
$H(G | \tilde{G}) < \infty$, and  $G$ and $\tilde{G}$ are mutually absolutely continuous.

\subsubsection{${\cal P}(\epsilon, \delta)$ and ${\cal P}_{a}(\epsilon)$ policies}

 Consider a ${\cal P}(\epsilon, \delta)$ policy acting
on samples  $(X_{i,j}, i=1,2, j \geq 1)$ from two distributions in  ${\cal L}$.

Suppose that under probability measure 
$P$, distribution of $(X_{1,j}, j \geq 1)$ is $F$,
and that of $(X_{2,j}, j \geq 1)$  is $G$.
Under  $\tilde{P}$,
the respective distributions are $F$ and $\tilde{G}$. 
Then,
\[
P(J=2) \leq \delta
\]
while
\[
\tilde{P}(J=2) \geq 1- \delta.
\]
Applying  (\ref{eqn:Kauff1}) to this set-up with  ${\cal E} = \{J=2\}$,
it follows that for $\delta \in (0, 1/2)$,

\begin{equation} \label{eqn:Kauff11}
E_P N_2(\tau)   \geq  \frac{d(\delta,1-\delta)}{H(G | \tilde{G})}.
\end{equation}
As noted by Kauffmann, Cappe and Garivier (2016)
\[
d(\delta,1-\delta) \geq \log \frac{1}{2.4 \delta}
\]
for all $\delta \in (0,1)$.

Theorem~\ref{theorem:thm_main}  is the main result of this section.  It requires 
 Lemma~\ref{lem:lem201} below.

\begin{lemma}  \label{lem:lem201}
Given  $G$ with finite mean $\mu_G$ and unbounded support on the positive real line, for any $\alpha >0$, and $\beta> \mu_G$
there exists a distribution $G_{\beta}$, such that the Kullback-Leibler distance between
 $G$ and $G_{\beta}$,
\begin{equation} \label{eqn:neg_result_odd_space1}
H(G | G_{\beta})   \leq \alpha,
\end{equation}
and
\begin{equation} \label{eqn:neg_result_odd_space2}
\mu_{G_{\beta}} = \int_{x \in \Re}x dG_{\beta}(x) \geq \beta.
\end{equation}
\end{lemma}

\begin{definition}
A collection of probability distributions ${\cal L}$  is referred to as $KL$ right dense, if each distribution has finite mean and  if for every $G\in {\cal L}$,
and every $\alpha>0$, $\beta> \mu_G$,  there exists a distribution  $G_{\beta} \in {\cal L}$
such that (\ref{eqn:neg_result_odd_space1}) and (\ref{eqn:neg_result_odd_space2}) hold.
\end{definition}

A necessary condition for ${\cal L}$ to be $KL$ right dense is that each member have an unbounded support
on the positive real line.
Examples of $KL$ right dense class includes collection of  distributions, which in addition to unbounded support
on the positive real line,  are also  light tailed
(recall that a distribution is said to be light tailed if its moment generating function is finite in a neighbourhood of zero).
Another example is  collection of distributions  that along with having an unbounded support
on the positive real line, are in ${\cal L}^p$, for some $p \geq 1$.  That is,  for some $p \geq 1$, their absolute $p$ th moment is finite.

\begin{theorem} \label{theorem:thm_main}
Under  ${\cal P}(\epsilon, \delta)$ policy operating on $KL$ right dense ${\cal {L}}$, for $P$ as above with distributions $F$ and $G$
in ${\cal L}$, for $\delta \in (0,1/2)$, 
\begin{equation} \label{eqn:lemma011}
E_a N_2(\tau)=  \infty. 
\end{equation}
Under ${\cal P}_{a}(\epsilon)$ policy operating on $KL$ right dense ${\cal {L}}$, for $P$ as above with distributions $F$ and $G$
in ${\cal L}$,
\begin{equation} \label{eqn:Kauff4}
\liminf_{\delta \rightarrow 0} \frac{E_P N_2(\tau)} {\log(1/\delta)} = \infty.
\end{equation}
\end{theorem}

\noindent {\bf Proof of Theorem \ref{theorem:thm_main}:}
Equation (\ref{eqn:lemma011}) follows from (\ref{eqn:Kauff1}) and 
  Lemma~\ref{lem:lem201}.

To see (\ref{eqn:Kauff4}), observe that  under ${\cal P}_{a}(\epsilon)$ policy,
\begin{equation} \label{eqn:Kauff001}
\limsup_{\delta \rightarrow 0} P(J=2) \delta^{-1} \leq 1,
\end{equation}
and
\begin{equation} \label{eqn:Kauff002}
\limsup_{\delta \rightarrow 0} (1-\tilde{P}(J=2)) \delta^{-1} \leq 1.
\end{equation}
Again applying  (\ref{eqn:Kauff1}) to this set-up with  ${\cal E} = \{J=2\}$,
it follows that 
\begin{equation} \label{eqn:Kauff3}
E_P N_2(\tau)   \geq  \frac{d(P(J=2)|\tilde{P}(J=2))}{H(G | \tilde{G})}.
\end{equation}

It is easy to see from (\ref{eqn:Kauff001}) and  (\ref{eqn:Kauff002}) that 
\[
d(P(J=2)|\tilde{P}(J=2)) \sim \log \delta^{-1}
\]
Hence,
\[
\liminf_{\delta \rightarrow 0} \frac{E_P N_2(\tau)} {\log(1/\delta)} \geq \frac{1}{H(G | \tilde{G})}.
\]
(\ref{eqn:Kauff4}) follows from Lemma~\ref{lem:lem201}.
$\Box$

\vspace{0.1in}

\begin{remark} \label{remark_on_unbounded_support} {\em
Observe that under ${\cal P}(\epsilon, \delta)$ policy, a bound on $E_P N_1(\tau)$  can be obtained by essentially repeating the analysis:
\begin{equation} \label{eqn:lemma0101}
 \frac{E_P N_1(\tau)} {\log(1/2.4 \delta)}
\geq  \sup_{\tilde{F} \in {\cal L}, \mu_{\tilde{F}} < \mu_G - \epsilon, H(F| \tilde{F})< \infty} \,\,\,
\frac{1}{H(F| \tilde{F})}.
\end{equation}
If  the objective of the proposed policies was to find a population with the smallest mean, one simply alters the sign
of the random variables in the above analysis. Then,  ${\cal L}$  needs to be $KL$ left dense (defined analogously to KL right dense collection) to ensure that
under the associated ${\cal P}(\epsilon, \delta)$ expected number of samples is infinite, while a 
${\cal P}_{a}(\epsilon)$ policy requires more than $O(\log (1/\delta))$ samples.
}
\end{remark}

\subsection{More general setting}
It is easy to extend this analysis in a variety of ways. For instance, suppose that given two populations our aim is to identify the one with
the largest  $p$ th quantile. Recall that for any real valued random variable $Z$ with distribution function $F_Z(\cdot)$, the quantile function is the inverse of $F_Z(\cdot)$, and in particular,  $p$ th quantile is given by
\[
F_Z^{-1}(p)= \inf_{x \in \Re} \{ F_Z(x) \geq p \}.
\]
The analysis easily generalises to handle such cases. To see this, consider random elements taking values
in a general state space ${\cal X}$.
Let $h$ denote a mapping from a probability distribution (more generally, a probability measure) on ${\cal X}$   to the real line.  The $p$ th quantile
is one example of such a mapping  from a probability distribution of a real-valued random variable to the real line.
Let ${\cal H}$ denote  a collection of probability distributions on  ${\cal X}$ such that $h(G)< \infty$ for all
$G \in {\cal H}$.  Suppose that in comparing two distributions $F$ and $G \in {\cal H}$, our aim is to identify $\max(h(F), h(G))$.

Now define   a  policy ${\cal P}(\epsilon, \delta)$ operating on distributions  in ${\cal H}$  with the property that
for   $F, G \in {\cal H}$ such that $|h(F)-h(G)| > \epsilon>0$ the policy ${\cal P}(\epsilon, \delta)$,  generates samples from the two distributions, adaptively deciding which distribution to sample from next and at some stage selects one of the two distributions as the one with the higher $h$ value, where $h$ may be evaluated at the respective empirical distribution.  The policy guarantees that  $P(FS) \leq \delta$.  As before, ${\cal P}_{a}(\epsilon)$ guarantees this asymptotically.

In the discussion at the beginning of Section~\ref{secs:anal}   where  the two populations  $(X_{i,j}: i = 1,2, j \geq 1)$ are  compared under the two probability measures
$P_a$ and $P_b$, replace $\mu_F, \mu_G$ and  $\mu_{\tilde{G}}$ with
$h(F), h(G)$ and $h(\tilde{G})$, where $F, G, \tilde{G} \in {\cal H}$. Then, it follows that under  ${\cal P}(\epsilon, \delta)$,
\begin{equation} \label{eqn:lemma01078}
 \frac{E_a N_2(\tau)}{\log(1/ 2.4\delta)}
\geq \sup_{ \tilde{G} \in {\cal H},
h(\tilde{G}) > h(F) + \epsilon, H(G|\tilde{G}) < \infty}\frac{1}{H(G| \tilde{G}) }.
\end{equation}
Under  ${\cal P}_{a}(\epsilon)$ policy
\begin{equation} \label{eqn:lemma01079}
\lim \inf_{\delta \rightarrow 0} \frac{E_a N_2(\tau)}{\log(1/ \delta)}
\geq \sup_{ \tilde{G} \in {\cal H},
h(\tilde{G}) > h(F) + \epsilon, H(G|\tilde{G}) < \infty}\frac{1}{H(G| \tilde{G}) }.
\end{equation}

\begin{example}{\em
In the case where ${\cal X}= \Re$, $h(F)$ denotes the $p$ th quantile of $F$, it easy to see that
RHS on (\ref{eqn:lemma01078}) can be infinite. We show this when ${\cal H}$ includes mixture of Gaussian distributions.

For $\mu>0$,  and $0< \epsilon < p <1/2$, consider the pdfs
\[
g(x) = \sqrt { \frac{1}{2 \pi}} \left ( p \exp( -x^2/2)   + (1-p) \exp(-(x-\mu)^2/2) \right ),
\]
and
\[
g_{\epsilon}(x) = \sqrt { \frac{1}{2 \pi}}
\left ( (p -\epsilon) \exp( -x^2/2)   + (1-p+\epsilon) \exp(-(x-\mu)^2/2) \right ).
\]

Since,
\[
\sqrt{\frac{1}{2 \pi}} \int_{\mu/2}^{\infty} \exp(-x^2/2) dx
=\sqrt{\frac{1}{2 \pi}} \int_{-\infty}^{\mu/2} \exp(-(x- \mu)^2/2) dx
\]
and  $p < 1/2$,
 the $p$ th quantile of $g$, $m_g < \mu/2$. The $p$ th quantile of $g_{\epsilon}$, $m_{g_{\epsilon}} \geq m(\epsilon)$
where $m(\epsilon)$ denotes the point where the second Gaussian component has a total mass equal
to $\epsilon$, that is,
\[
 \frac{\epsilon}{1-p + \epsilon} =  \sqrt{\frac{1}{2 \pi}} \int_{-\infty}^{m(\epsilon)} \exp(-(x- \mu)^2/2) dx
 = \sqrt {\frac{1}{2 \pi}} \int_{-\infty}^{(m(\epsilon) -\mu)} \exp(-x^2/2) dx.
 \]
Using the well known fact that
\[
\int_{-\infty}^{z} \exp(-x^2/2) dx \sim \frac{1}{|z|}\exp(-z^2/2)
\]
as $ z \rightarrow -\infty$, we have
\[
\mu - m(\epsilon) \sim  \sqrt{2 \log ((1-p+\epsilon)/\epsilon)}
\]
as $\epsilon \rightarrow 0$.
In particular, it follows that for small and fixed $\epsilon$, $m_{g_{\epsilon}} - m_g$ can be made arbitrarily large by increasing $\mu$.

Now consider, the Kullback-Leibler distance
\begin{equation} \label{eqn:ent_once_more}
\int_{-\infty}^{\infty} \log \left (   \frac{g(x)}{g_{\epsilon}(x)} \right )g(x) dx.
\end{equation}
Observe that
\[
\left ( \frac{p}{p-\epsilon} \right ) g_{\epsilon}(x) \geq g(x),
\]
for all $x$,  for all $\mu$. Thus,
$\frac{g(x)}{g_{\epsilon}(x)} \leq \frac{p}{p-\epsilon}$ for all $x$ and $\mu$. Hence, for any $\mu$,  (\ref{eqn:ent_once_more})
can be made arbitrarily close to zero by choosing $\epsilon$ sufficiently close to zero.
}
\end{example}

\section{Positive results - the non-adaptive algorithms}

In this section, we show that under conditions on moments of strictly increasing non-negative  convex functions of the underlying populations,
for any given $\epsilon, \delta >0$, one can develop  non-adaptive  ${\cal P}(\epsilon, \delta)$ algorithms that require
a deterministic and known  $O(\log(1/\delta))$  computational effort for any instance of underlying populations. These algorithms rely on truncating or capping the random samples generated and carefully bounding the resulting worst-case bias using explicitly available moment bounds.

Specifically, we suppose (as in the Introduction) that we can sample independently from $d$ different random variables $(X(i): i \leq d)$. Further, from each population $i$, we can generate independent identically distributed samples $(X(i,j): j \geq 1)$.  Recall that the distribution of $(X(i): i \leq d)$   is unknown
and our aim is to find
\[
 i^* = \arg \max_{1 \leq j \leq d} EX(j).
\]

First in Section~\ref{subsec:01},  we assume that  $X(i) \in [0,b]$ a.s for each $i \leq d$, for some $b>0$, and review the
${\cal P}(\epsilon, \delta)$ algorithms  that require $O(\log(1/\delta))$  computational effort in that simple setting.
Analysis is straightforward and is discussed (independently, it appears), e.g.,  in Even-Dar et. al. (2002, 2006), Glynn and Juneja (2004, 2011).
In Section~\ref{subsec_bound_trunc}, we arrive at  the maximum expected error that may result through
either appropriately truncating or capping a random variable  when an upper bound on the moment of strictly increasing, non-negative  convex function of a random variable is explicitly known. In Section~\ref{subsec:02}, we develop
${\cal P}(\epsilon, \delta)$ algorithms that require $O(\log(1/\delta))$ computational effort,
when appropriate moment upper bounds are available for each population.  Such bounds can often be found in simulation models by the use of Lyapunov function based techniques (see, e.g., Glynn and Zeevi 2008).

\subsubsection{${\cal P}(\epsilon, \delta)$ policy for bounded random variables} \label{subsec:01}

For $\epsilon, b>0$, consider ${\cal X}_{\epsilon}(b) = \{(X(i): i \leq d): EX(i^*) > EX(j) +\epsilon \,\,\ \forall j \ne i^*,   X(i) \in [0, b] \,\,
\forall i \}$.

\vspace{0.1in}

A reasonable
algorithm on ${\cal X}_{\epsilon}(b)$ for a well chosen $n$ (discussed later) is:

\begin{itemize}
\item
Generate independent samples $(X(i,j): i=1, \ldots, d \mbox{ and } j = 1, \ldots,n)$.
\item
Let $\bar{X}(i) = \frac{1}{n} \sum_{j=1}^n X(i,j)$.  Declare
\[
 \hat{i} = \arg \max_{1 \leq i \leq d} \bar{X}(i)
\]
as the best design.
\end{itemize}

Recall that false selection occurs if $ \hat{i} \ne i^*$,  with probability
\[
P(\bar{X}(i^*) < \max_{j \ne i^*} \bar{X}(j)).
\]
This is bounded from above by
\[
\sum_{j \ne i^*} P(\bar{X}(i^*) < \bar{X}(j))
\]
Using Hoeffding's inequality, we have
\[
 P(\bar{X}(i^*) < \bar{X}(j)) \leq
 P(\bar{X}(i^*) - \bar{X}(j) - (EX(i^*) - EX(j)) < - \epsilon)
 \leq \exp \left (- n \frac{\epsilon^2}{2 b^2} \right ).
 \]
Thus, $n =\frac{2 b^2}{\epsilon^2} \log ((d-1)/\delta)$ provides the desired ${\cal P}(\epsilon, \delta)$ policy
for any set of distributions corresponding to  ${\cal X}_{\epsilon}(b)$.

\subsection{ Bounding the worst-case truncation or capping bias} \label{subsec_bound_trunc}

Suppose that
${\cal X} $ is a
class of non-negative random variables and $f$ is a strictly increasing non-negative convex function.
Examples include $f(x)= x^{\alpha}$ for $x \geq 0$ and $\alpha >1$, and $f(x) = \exp(\theta x)$ for $\theta >0$.
We discuss two formulations to bound worst-case biases resulting from 1) truncating a random variable, and 2) capping it.

Consider the optimization problem ${\bf O_1}$
\begin{eqnarray} \label{eqn:1}
 & \max_{X \in {\cal X} }   E X  I(X \geq u)   \\
\mbox{  such that }  & E f(X)   \leq c,  \label{eqn:2}
\end{eqnarray}
for some positive $u$ and $c$.

Also consider ${\bf O_2}$, where the objective function is
instead set to
\begin{equation}
\max_{X \in {\cal X} }   E (X -u)I(X > u),
\end{equation}
again under constraint (\ref{eqn:2}).

Furthermore, since,  by Jensen's
inequality,
 \[
 E f(X) \geq f(E X ) \geq f(0),
 \]
 we assume that  $c > f(0)$.

 ${\bf O_1}$ denotes the worst-case expected truncation bias under constraint (\ref{eqn:2})
 when $X$ is replaced by the truncated $XI(X < u)$.  ${\bf O_2}$ denotes the smaller bias
 when $X$ is replaced by the capped $\min(X,u)$.

We now argue that the solutions to ${\bf O_1}$ and ${\bf O_2}$
put their mass at most at two points. To see this,  observe that given any non-negative random variable $X$ that satisfies
(\ref{eqn:2}), a two-valued random variable $Y$
that takes value
\begin{itemize}
\item
$E[X| X < u]$ ($E[X| X \leq u]$) with probability
$P(X < u)$ ($P(X \leq u)$),  and
\item
value  $E[X| X \geq u]$ ($E[X| X > u]$) with probability
 $P(X \geq u)$ ($P(X > u)$),
\end{itemize}
has  the same mean $EY=EX$ and same objective function
values under ${\bf O_1}$ (${\bf O_2}$) , i.e.,
\[
  EYI(Y \geq u) = EXI(X \geq u), \,\,\,\,\,
\left (E(Y-u)I(Y > u) = E(X-u)I(X > u) \right ).
\]
Furthermore,  $Ef(Y) \leq Ef(X)$, with equality only if
$X=Y$ a.s.  Thus, only random variables that take at most two values  can solve
our optimization problems ${\bf O_1}$ and ${\bf O_2}$. It is also easy to see that at the optimal solution in both the cases, the constraint
(\ref{eqn:2}) has to be tight.  Hence, we restrict our search to
random variables
 taking values $0 \leq x_1 < x_2$  with probability
$1-p$ and $p  \in [0,1]$ where
\[
(1-p)f(x_1)  + p f(x_2) = c,
\]
 so that $0 \leq x_1 < f^{-1}(c)$ and $ x_2 \geq f^{-1}(c)$. Furthermore,
 \[
 p = \frac{c- f(x_1)}{f(x_2) - f(x_1)}.
 \]

Propositions~\ref{prop:july1115-1} and \ref{prop:July1115} below observe that
 unique (a.s.) solutions $\tilde{X}_1$  and $\tilde{X}_2$, respectively,  to the optimization problems
 ${\bf O_1}$ and ${\bf O_2}$ are either degenerate and equal $ f^{-1}(c)$ with probability 1, or they take value zero with positive probability.

\begin{proposition} \label{prop:july1115-1}
The unique optimal solution $\tilde{X}_1$ for ${\bf O1}$,
\begin{enumerate}
\item
equals $f^{-1}(c)$ with probability 1, for $u \leq f^{-1}(c)$.
\item
For $u > f^{-1}(c)$,
$\tilde{X}_1$ has a two-value distribution. It equals $u$ with probability
\[
\frac{c-f(0)}{f(u) -f(0)},
\]
and zero otherwise.
\end{enumerate}
\end{proposition}

The following assumption considerably eases the analysis of ${\bf O2}$.

\begin{assumption} \label{ass_july1115}
Given any $u>0$, there exists
$x_u>u$, the unique  solution to
\[
x-u = \frac{f(x) - f(0)}{f'(x)}.
\]
\end{assumption}

The above assumption can be seen to hold , e.g., for $f(x) =  x^{\alpha}$, $\alpha >1$ and
for $f(x) = \exp(\theta x)$, $\theta >0$. To see this, observe that
in these cases the function
$f(x)-f(0) - (x-u)f'(x)$ monotonically decreases from a positive value at $u$ to $-\infty$ as
$x \uparrow \infty$.

\begin{proposition} \label{prop:July1115}
The unique optimal solution $\tilde{X}_2$ for ${\bf O2}$ under Assumption~\ref{ass_july1115},
\begin{enumerate}
\item
a.s. equals  $ f^{-1}(c)$ for $x_u \leq f^{-1}(c)$.
\item
For $x_u > f^{-1}(c)$,
$\tilde{X}_2$ has a two-valued distribution. It equals $x_u$ with probability
\[
\frac{c-f(0)}{f(x_u) -f(0)},
\]
and zero otherwise.
\end{enumerate}
\end{proposition}

\begin{remark} \label{remark:rem1} {\em

Suppose that  $f(x) = x^{\alpha}$, $\alpha > 1$. Under ${\bf O_1}$,  the solution corresponds to
 $\tilde{X}_1=c^{1/\alpha}$ with probability 1 for $u \leq c^{1/\alpha}$ and the associated objective function value
 equals $ c^{1/\alpha}$. Otherwise,
$\tilde{X}_1$ takes two values 0 and $u$ where the probability of the latter equals $c u^{-\alpha}.$
The objective function value then equals $c u^{-(\alpha-1)}$.

Under ${\bf O_2}$, for $u>0$
\[
x_u = u \left ( \frac{\alpha}{\alpha-1} \right ).
\]
Then,  $\tilde{X}_2 =c^{1/\alpha}$ with probability 1 for $u \leq c^{1/\alpha}\left ( \frac{\alpha-1}{\alpha} \right )$,
and the associated objective function value equals
 $c^{1/\alpha}-u$. Otherwise,
$\tilde{X}_2$ takes two values  0 and $x_u$ where the probability of the latter equals
\[
 c u^{-\alpha}\left ( \frac{\alpha-1}{\alpha} \right )^{\alpha}.
\]
The optimal objective function value equals
\begin{equation}  \label{eqn:error_change}
  c u^{-(\alpha-1)}\left ( \frac{(\alpha-1)^{\alpha-1}  }{\alpha^{\alpha}} \right ).
\end{equation}
Thus the worst case bias reduces by factor  $\left ( \frac{(\alpha-1)^{\alpha-1}  }{\alpha^{\alpha}} \right )$   if we use $\min(X,u)$ instead of $XI(X <u)$
to get random variables bounded from above by $u$.

Also note that under ${\bf O_1}$ if
$f(x) = \exp(\theta x)$
for some $\theta >0$, then  $\tilde{X}_1= \log c/\theta$ for
$u \leq \log c/\theta$ with probability 1, and
$\tilde{X}_1$ takes two values  0 and $u$ where the probability of the latter equals
\[
  \frac{c- 1}{\exp(\theta u)-1}.
\]
The objective function value in this case equals $u  \frac{c- 1}{\exp(\theta u)-1}$.
 }
\end{remark}
\subsubsection{Threshold function}

Going forward, we focus on using capping to bound the random variables, because this has smaller bias compared to truncating.
Let $g(u)$ denote the optimal value of ${\bf O2}$ as a function of $u$.
Since,
\[
E (X -u)I(X > u) = \int_{u}^{\infty} P(X>y) dy,
\]
$g(u)$ is a strictly non-increasing  function of $u$ when $X$ is unbounded. Let the {\em threshold function} $r(\cdot)$ denote the inverse of
$g(\cdot)$. Thus, if wish to restrict the bias to $\epsilon$, let the { threshold function} $r(\epsilon)$ denote the
value of $u$ in ${\bf O2}$ for which the optimal value $g(u)$  equals $\epsilon$.

 \subsection{Analysis for unbounded random variables} \label{subsec:02}

We now return to the problem of finding ${\cal P}(\epsilon, \delta)$ algorithms that require
a deterministic and known  $O(\log(1/\delta))$  computational effort for any instance of underlying populations, when
upper bounds on the moments of strictly increasing non-negative convex functions of the underlying random variables  are explicitly known.
Specifically, let ${\cal Y}$ denote a collection probability distributions corresponding to non-negative random variables $(X(i): i \leq d)$
such that there exists strictly convex, non-decreasing functions
$(f_i(\cdot): i \leq d)$ and positive constants  $(c_i: i \leq d)$ so that
\[
E f_i(X(i)) \leq c_i
\]
for all $i$. Further, let $r_i(\cdot)$ denote the  threshold functions associated with each $f_i(\cdot)$.
Set
\[
R(x) = \max_{i \leq d} r_i(x).
\]

Consider the following algorithm for selecting the population with the highest mean.
Let $\beta \in (0,1)$ and $ \epsilon >0$.
\begin{itemize}
\item
Generate independent samples $(X(i,j): i=1, \ldots, d \mbox{ and } j = 1, \ldots,n_{\beta})$.
Let $Y(i,j) = \min(X(i,j), R(\beta \epsilon))$ for  all $i,j$.
\item
Let $\bar{Y}(i) = \frac{1}{n_{\beta}} \sum_{j=1}^{n_{\beta}} Y(i,j)$.  Declare
\[
 \hat{i} = \arg \max_{1 \leq i \leq d} \bar{Y}(i)
\]
as the best design.
\end{itemize}

\begin{proposition} \label{prop:1what ye lose}
The above algorithm is   ${\cal P}(\epsilon, \delta)$
for distributions in ${\cal Y}$ with
\[
n_{\beta} =  \lceil \frac{2 R(\beta \epsilon)^2}{\epsilon^2 (1-\beta)^2} \log ((d-1)/\delta)  \rceil
\]
providing the desired guarantees.
\end{proposition}

{\noindent \bf Proof of Proposition~\ref{prop:1what ye lose}:}
The probability of false selection corresponds to
\begin{equation} \label{eqn:prob_false_selection_pos}
P(\bar{Y}(i^*) < \max_{j \ne i^*} \bar{Y}(j)).
\end{equation}
Repeating the analysis in Section~\ref{subsec:01}, keeping in mind that
for $ j \ne i^*$,
\begin{eqnarray*}
E {Y}(i^*) - E{Y}(j)
& =  & E \min(X(i^*), R(\beta \epsilon)) -
E \min(X(j), R(\beta \epsilon))  \\
&  \geq &  E X(i^*) - \beta \epsilon  - E X(j)  \\
& \geq  & (1-\beta) \epsilon,
\end{eqnarray*}
we conclude that (\ref{eqn:prob_false_selection_pos}) is bounded from above
by $\delta$.
$\Box$

\begin{remark}
{\em
The  $\beta$ that minimises $n_{\beta}$ corresponds to that minimising
$\frac{R(\beta \epsilon)^2}{1-\beta)^2}$.

We solve this when  $f_i(x) = x ^{\alpha}$, $\alpha >1$ for all $i$.
Then,  it can be seen from Remark~\ref{remark:rem1}
 that
 \[
 r_i(x) =   \left ( \frac{c_i}{x} \right )^{1/(\alpha-1)}\left ( \frac{\alpha-1}{\alpha^{\alpha/(\alpha-1)}} \right )
 \]
 and we can set
 \[
 R(x) = \frac{1}{x^{1/(\alpha-1)}}\left ( \frac{\alpha-1}{\alpha^{\alpha/(\alpha-1)}} \right )  \max_{i \leq d} {c_i}^{1/(\alpha-1)}.
\]
Then minimising $n_{\beta}$ corresponds to maximising
\[
\beta^{2/(\alpha-1)}(1-\beta)^2,
\]
and is achieved at $\beta =1/\alpha$.
Thus the bound on the number of samples needed equals
$n_{1/\alpha}$ or
\[
\frac{2 R(\beta \epsilon)^2}{\epsilon^2 (1-\beta)^2} \log ((d-1)/\delta)
=2 \left ( \frac{  \max_{i \leq d}c_i  }{\epsilon^{\alpha}}\right )^{2/(\alpha-1)}
\log ((d-1)/\delta).
\]
}
\end{remark}

\section{Sequential pure exploration algorithm}

We now review some of the related literature that comes under the broad topic of stochastic multi-arm bandit
methodology. We discuss the elegant sequential sampling strategy referred to as the successive elimination algorithm proposed by Even-Dar et. al. (2002, 2006).
Although they also propose a slightly more effective median elimination algorithm in that paper, and there have been significant developments on the pure exploration problem since then (see, e.g., Audibert and Bubeck 2010, Jamieson  et. al. 2013, Garivier and Kaufmann 2016), the sequential algorithm of Even-Dar et. al. (2006) is particularly simple and lends to an easier laconic discussion.  They considered the setting where the underlying random variables were Bernoulli, while we allow generally distributed random variables when explicit  bounds on the moments are available. As mentioned in the introduction, we use one of the proposed methods  from Bubeck et. al. (2013) for this purpose that  relies on careful truncation of underlying random variables. Their analysis focusses on the regret minimisation objective but is easily adapted to our pure exploration setting. We present a minor tweak - while they considered truncations of the form $XI( X <u)$ to bound rv $X$, we  note the minor performance benefits from using instead the capped random variable $\min(X,u)$. These may also extend to the regret minimisation objective. In addition, we compute bounds on the expected number of samples generated under the pure exploration algorithm for general random variables.

\subsection{Sequential algorithm}

We refer to populations as arms in this section.
Thus, we have $d$ arms. Arm $i$ when pulled gives a reward distributed as $X(i)$ and our aim is to find the best arm
$i^* = \max_{i \leq d} EX(i).$
Let  the maximum be achieved by a unique arm and let $\Delta_i = EX(i^*) - EX(i) >0$
for all $i \ne i^*$.

Suppose that there exists a non-negative function $\alpha_m$ with the property that for $i \leq d$,
\begin{equation} \label{eqn:use_hoeffding}
P(| \bar{X}_m(i) - EX(i) | > \alpha_m) \leq \frac{c \delta}{m^2 d}
\end{equation}
where $\bar{X}_m(i)$ denotes the average of $m$ i.i.d. samples of $X(i)$,
and
\[
c = \left (\sum_{m=1}^{\infty}1/m^2 \right )^{-1} = 6/\pi^2.
\]
In the examples that we consider later $\alpha_m$ is seen to be decreasing for
all $m \geq e$. In our discussions, we will ignore this minor issue and assume that
 $\alpha_m$ is a decreasing, non-negative function of $m$.

Equation (\ref{eqn:use_hoeffding})  ensures that
\[
P(E_{\delta}) \geq 1- \delta,
\]
where
\[
E_{\delta}= \{| \bar{X}_m(i) - EX(i) | \leq \alpha_m, \forall m, \, \forall i \leq d \},\]
and is the rationale for the successive elimination algorithm outlined below.

\vspace{0.1in}

{\noindent \bf Successive elimination algorithm}

\begin{enumerate}
\item
Set $m=1$ and $S=\{1, 2 \ldots, d\}$.
\item
Set for each arm $i$, $\bar{X}_1(i)=0;$
\item
{\bf Repeat}
 \begin{itemize}
\item
Sample every arm $i \in S$ once and let $\bar{X}_m(i)$ be the average reward of arm $i$
by trials or pulls $m$;
 \item
 Let $\bar{X}_m(\max)= \max_{i \in S}\bar{X}_m(i)$;
 \item \hspace{0.5in}
 For each arm $i \in S$ such that $\bar{X}_m(\max) - \bar{X}_m(i) \geq 2 \alpha_m$ {\bf do}
 \item \hspace{1in} set $S= S-\{i\}$;
 \item \hspace{0.5in} end
 \item $m =  m+1$;
  \end{itemize}
{\bf Until $|S| >1$};
\end{enumerate}

It is easy to see that on the set $E_{\delta}$ the best arm is never eliminated and that all other arms
are eventually eliminated (see Even-Dar et. al. 2002, 2006).  Also, It can be easily checked
using Hoeffding's inequality  that $\alpha_m= b\sqrt{\frac{2}{m} \log \left (\frac{d m^2}{c \delta} \right )}$ works for random variables $X(i) - EX(i) \in [-b, +b]$, $i \leq d$  in (\ref{eqn:use_hoeffding}).

\subsubsection{Explicit moment bounds on random variables}

The key step above was simply to arrive at the set
 $E_{i, \delta}$
\[
= \{|\bar{X}_m(i) - EX(i)| < \alpha_m \mbox{ for all m}\},
\]
that has probability at least  $1-\delta/d$.

Bubeck et. al. (2013)  achieve this when each $X(i)$
has an explicit bound on its moment. They focus on the on the more interesting case
where the bound is  on $\alpha$ moment for $\alpha \in (1,2]$. We also restrict our discussion
to this case.  Suppose that $K$ upper bounds an $\alpha$ abolute moment of a random variable for $\alpha >2$. Then,
$K^{\alpha/2}$ upper bounds its second moment. Similarly, an upper bound on an exponential moment
of a random variable can be reduced to an upper bound on the second moment. 

A key to   analysis of Bubeck et. al. (2013) is Lemma~\ref{Lemma:1} below that relies on using truncation and Bernstein inequality (shown with the proof of   Lemma~\ref{Lemma:1}). In Lemma~\ref{Lemma:1}, we also state the result when the random variables are capped instead of truncated. For this case, we assume that the underlying random variables are non-negative.
For $\alpha >1$, let $\tilde{\alpha} = \left (\frac{(\alpha-1)}{\alpha} \right )^{\alpha-1}$.

\begin{lemma} \label{Lemma:1}
Let $\delta \in (0,1)$, $\alpha \in (1,2]$, $K>0$,
and $(X_i: i \leq n)$ be iid samples of rv $X$.
Suppose that $E |X|^{\alpha} \leq K$ and let
\[
B_m = \left (  \frac{K m}{\log (\delta^{-1})}  \right )^{\alpha^{-1}}.
\]
Then, with probability
at least $1-\delta$,
\begin{equation} \label{eqn_pos_part_1}
\bar{X}_n \geq EX  - (2+\alpha)  K^{ 1/\alpha}
\left (  \frac{\log (\delta^{-1})}{n}   \right )^{\frac{\alpha-1}{\alpha}},
\end{equation}
where
\[
\bar{X}_n
= \frac{1}{n} \sum_{m=1}^n X_m I(|X_m| \leq B_m)
\]
denotes the empirical mean of truncated samples.

When $X$ is also assumed to be non-negative, then again with probability
at least $1-\delta$,
\begin{equation} \label{eqn_pos_part_2}
\tilde{X}_n \geq EX  - (2+\tilde{\alpha})  K^{ 1/\alpha }
\left (  \frac{\log (\delta^{-1})}{n}   \right )^{\frac{\alpha-1}{\alpha}},
\end{equation}
where
\[
\tilde{X}_n
= \frac{1}{n} \sum_{m=1}^n \min( X_m, B_m),
\]
denotes the capped empirical mean.
\end{lemma}

The reverse of (\ref{eqn_pos_part_1})
\[
\bar{X}_n \leq EX  + (2+\alpha)  K^{ (1+\epsilon)^{-1} }
\left (  \frac{\log (\delta^{-1})}{n}   \right )^{\frac{\alpha-1}{\alpha}}.
\]
also holds with probability at least $1-\delta$ following essentially identical proof. Similarly, the reverse of
(\ref{eqn_pos_part_2}).

\subsection{Expected number of samples generated}

Lemma below is useful to our analysis.

\begin{lemma} \label{lemma:1123}
Suppose that
$a \geq e, b \geq 1$ and $t=t^* \geq 1$ solves
\[
a+ b \log t = t.
\]
Then,
\[
t^* \leq a + b \log a + \frac{2b^2}{a} \log(a+b).
\]
\end{lemma}

To compute $E T(i)$ the expected number of times  arm $i \ne i^*$ is pulled, note that
\[
E T(i) = \sum_{m=1}^{\infty} P( \bar{X}_m(\max) - \bar{X}_m(i) < 2 \alpha_m)
\leq \sum_{m=1}^{\infty} P( \bar{X}_m(i^*) - \bar{X}_m(i) < 2 \alpha_m).
\]
For $i \ne i^*$, let
\[
\tau^*_i = \inf \{ m: 4 \alpha_m \leq  \Delta_i \}
\]
(recall that $\Delta_i = E X(i^*) - EX(i)$).
Then,
\begin{eqnarray*}
E T(i)  & \leq & \tau^*_i + \sum_{m= \tau^*_i+1}^{\infty}
P( \bar{X}_m(i^*) - EX(i^*)  - (\bar{X}_m(i) - EX(i))
 < - 2 \alpha_m), \\
 & \leq & \tau^*_i + \sum_{m=\tau^*_i+1}^{\infty} \left ( P(\bar{X}_m(i^*) \leq EX(i^*) -\alpha_m)
 +   P(\bar{X}_m(i) \geq EX(i) +\alpha_m)  \right ) \\
 & \leq & \tau^*_i + 2 \delta/d.
\end{eqnarray*}

Thus, the total expected number of samples generated
for all arms $i \ne i^*$ is bounded from above by
\[
\sum_{i \ne i^*} \tau_i  + 2 \delta
\]
and the total number of samples generated
is bounded from above by twice this amount.

\subsubsection{ Bounded random variables}

In particular, when
$X(i) - EX(i) \in [-b, +b]$,  $\alpha_m= b\sqrt{\frac{2}{m} \log \left (\frac{d m^2}{c \delta} \right )}$ works.
Let $m^*$ be the solution to
\[
4b\sqrt{\frac{2}{m} \log \left (\frac{d m^2}{c \delta} \right )} =  \Delta_i.
\]
Then, $\tau^*_i \leq m^*+1$.

Hence, using Lemma~\ref{lemma:1123},
\[
\tau^*_i \leq \frac{32b^2}{\Delta_i^2} \log \left (\frac{d}{c \delta} \right )
 + \frac{64b^2}{\Delta_i^2} \log \left (    \frac{32b^2}{\Delta_i^2} \log \left (\frac{d}{c \delta} \right )   \right )
 +1
 \]
plus terms that become small as $\delta$ decreases to zero.

Recall that the total number of samples generated is bounded from above by
\[
2 \sum_{i \ne i^*} \tau^*_i  + 4 \delta.
\]
Hence the dominant terms in the upper bound for total number of expected samples for small $\delta$ are
\[
64 b^2 \log \left (\frac{d}{c \delta} \right ) \sum_{i \ne i^*}\frac{1}{\Delta_i^2} .
\]

\subsubsection{ Explicit bound on moments}

When $E|X(i)|^{\alpha} \leq  K$ for all $i \leq d$, it is easily seen that
\[
\alpha_m = (2+\tilde{\alpha}) K^{1/\alpha} \left ( \frac{\log \left ( \frac{2m^2 d}{c \delta}\right )}{m}  \right )^{(\alpha-1)/\alpha}
\]
satisfies (\ref{eqn:use_hoeffding}).

Again, the number of times all the arms  $i \ne i^*$ are pulled is bounded from above by
\[
\sum_{i \ne i^*} \tau^*_i  + 2 \delta
\]
where now $\tau^*_i$ is bounded from above by $m^*+1$ and $m^*$ solves the equation
\[
4 (2+\tilde{\alpha}) K^{1/\alpha} \left ( \frac{\log \left ( \frac{2m^2 d}{c \delta}\right )}{m}  \right )^{(\alpha-1)/\alpha}
= \Delta_i.
\]

From Lemma~\ref{lemma:1123}, it follows that
\[
m^* \leq a + b \log a + \frac{2b^2}{a} \log(a+b)
\]
where
\[
a =  \left (\frac{4 (2+\tilde{\alpha}) K^{1/\alpha}}{\Delta_i} \right )^{\frac{\alpha}{\alpha-1}} \log \left ( \frac{2d}{c \delta}\right )
\]
and
\[
b = 2 \left (\frac{4 (2+\tilde{\alpha}) K^{1/\alpha}}{\Delta_i} \right )^{\frac{\alpha}{\alpha-1}}.
\]
Hence, the dominant terms in the upper bound for total number of expected samples for small $\delta$ are
\[
2 \left (4(2+\tilde{\alpha}) K^{1/\alpha} \right )^{\frac{\alpha}{\alpha-1}} \log \left ( \frac{2d}{c \delta}\right )
   \sum_{i \ne i^*} \left ( \frac{1}{\Delta_i} \right )^{\frac{\alpha}{\alpha-1}} .
\]

\section{Conclusions}

We considered the well studied classical problem of identifying the best population  amongst many when the populations are ranked by the expectations of associated random variables. These expectations may not be known, but samples of the underlying random variables can be generated, e.g.,  via simulation.  In the existing literature, when the underlying random variables are light-tailed, the probability of false selection is known to decay at an exponential rate. The question that we addressed in the paper was
whether this convergence rate could be exploited to build fast algorithms  with desired probabilistic guarantees.
Specifically, can algorithms that require
$O(\log (1/\delta))$ computational effort be constructed that  restrict the probability of false selection to a small $\delta$.
The existing literature had suggested  that by empirically estimating large deviations rate function
of underlying random variables, one may be able to  develop such algorithms.  
We  answered this question in negative when the underlying random variables belong to a class of distributions whose support on the right is unbounded and that are right KL-dense. 

The key conclusion
from our analysis  is that further restrictions are needed to develop algorithms that require
$O(\log (1/\delta))$ computational effort and provide explicit  or asymptotic guarantees of restricting the probability of false selection to a small $\delta$. We showed that explicit bounds on expectations of convex functions of underlying random variables suffices.
We also showed that sequential methods developed in multi-armed bandit literature for regret minimization can be adapted
to our pure exploration settings once moment bounds on underlying random variables are available.

\section{Appendix: Proofs}

{\noindent {\bf Proof of Lemma~\ref{lem:lem201}:}}
Consider a large $b$ whose value will be fixed later.
Furthermore, take $\gamma \in (0,1)$.
Construct  a probability distribution $G_k$ as follows: Set
\[
G_k(x) = (1-\gamma) G(x)
\]
for all $x \leq b$, and,
\[
\bar{G}_k(x) = \beta \bar{G}(x)
\]
for $ x >b$, where,
$\bar{G}_k(x)= 1- G_k(x)$ and
$\bar{G}(x)= 1- G(x)$.

Note that
\[
\beta = 1+ \gamma \frac{G(b)}{\bar{G}(b)} > 1.
\]

Then,
\begin{equation} \label{eqn:101}
0 \leq \int_{x \in \Re} \log \left (\frac{dG(x)}{dG_k(x)} \right ) dG(x)
\leq   - G(b) \log (1-\gamma) .
\end{equation}

By selecting
$\gamma = 1-\exp(-\alpha)$, we get
\[
- G(b) \log (1-\gamma) \leq \alpha.
\]
%Furthermore, then $\bar{G}(b) \log \beta$ equals
%\[
%\bar{G}(b) \log \left ( \frac{1- \exp(-\alpha/2) G(b)}{\bar{G}(b)} \right )
%\leq  \bar{G}(b)\log \left (1+ \frac{\alpha G(b)}{2 \bar{G}(b)} \right )
%< \alpha/2.
%\]
%(Using $e^{-x} \geq 1-x$ and $\log (1+x) \leq x$ above).

Also, for $b$ such that
$G(b^+) = G(b^-)$,
\[
\mu_{G_k} = (1-\gamma) \int_{-\infty}^{b} x d G(x)
+\left (1+ \gamma \frac{G(b)}{\bar{G}(b)} \right ) \int_{b}^{\infty} x d G(x)
\geq
\exp(-\alpha) \mu_G + (1- \exp(-\alpha)) G(b) b
\]
Since, RHS increases to infinity as
$b \rightarrow \infty$, one can select $b$ sufficiently large so that
$\mu_{G_k} \geq k$.
$\Box$

\vspace{0.1in}

{\noindent \bf Proof of Proposition~\ref{prop:july1115-1}:}
 The case $u \leq f^{-1}(c)$,
 follows by noting that
 \[
 f(EX) \leq E f(X) \leq c,
 \]
 so that
 $E X I(X \geq u) \leq EX \leq f^{-1}(c).$

Now consider
$u > f^{-1}(c)$. Suppose that $\tilde{X}_1$
takes values $ x^*_1 \in [0, f^{-1}(c))$ and
 $x_2^* \geq u$. ($x_2^* < u$ is not possible as then the optimal objective function value is zero.)

The objective function (\ref{eqn:1}) at optimal value then equals
\[
x_2^* \left (\frac{c- f(x_1^*)}{f(x_2^*) - f(x_1^*)} \right ).
\]

Note that
\[
  \frac{c - f(x)}{f(x_2^*)-f(x)}
\]
for $ x < f^{-1}(c)$ and  $f(x_2^*) > c$,
is a strictly decreasing function of $x$. Thus,
$x^*_1=0$.

Now observing that
\[
f(y) - f(0) < y f'(y)
\]
from simple calculus, it follows that
\[
y \left (\frac{c- f(0)}{f(y) - f(0)} \right )
\]
is non-increasing in $y$ for $y \geq u$ and
the result follows.
$\Box$

\vspace{0.1in}

{\noindent \bf Proof of Proposition~\ref{prop:July1115}:}
 Consider a solution
 that puts mass at two points
 $x^*_2$ and $ x^*_1 \in [0, x^*_2)$ with probabilities
 $p$ and $1-p$, respectively.
Clearly, then $x^*_1 < f^{-1}(c)$ and $x^*_2 \geq f^{-1}(c)$.

 First suppose that an optimal solution has $x^*_1 \geq u$ with positive probability.
 We show that this leads to a contradiction. Note that $x^*_1 \geq u$ only if
$u < f^{-1}(c)$.
In that case, the objective function equals
\begin{equation} \label{eqn:july1115}
\left ((x^*_2 -x^*_1)
 \frac{c - f(x^*_1)}{f(x^*_2)-f(x^*_1)}  + (x^*_1 -u) \right )
\end{equation}

Consider the function
\[
\left ((x -x^*_1)
 \frac{c - f(x^*_1)}{f(x)-f(x^*_1)}   \right ).
\]
This is clearly non-increasing in $x$. Thus $x^*_2$ in (\ref{eqn:july1115}) equals
$f^{-1}(c)$ with probability 1 providing the desired contradiction.

Now suppose that
$x^*_1 < u$. Then, to achieve a positive value of the objective function, we need
 $x^*_2 > u$. In particular, the optimal objective function equals
 \begin{equation} \label{eqn:july11152}
\left ((x^*_2 -u)
 \frac{c - f(x^*_1)}{f(x^*_2)-f(x^*_1)} \right ).
\end{equation}

Note that
\[
  \frac{c - f(x)}{f(x^*_2)-f(x)}
\]
for $ x < \min(x^*_2, f^{-1}(c))$ and  $f(x^*_2) \geq c$,
is a strictly decreasing function of $x$. Thus,
$x^*_1=0$ in  (\ref{eqn:july11152}).

Now consider the function
\[
\left ((x -u)
 \frac{c - f(0)}{f(x)-f(0)} \right )
\]
for $x > u$. Due to Assumption~\ref{ass_july1115}, this is maximised at $x_u$.
Thus, if $x_u > f^{-1}(c)$, then clearly $x^*_2 =x_u$ in (\ref{eqn:july11152}).
Else, $x^*_2 = f^{-1}(c)$ with probability 1. The proposition stands proved.
$\Box$

\vspace{0.2in}

The Bernstein inequality below is useful for the proof of Lemma~\ref{Lemma:1}.

\begin{lemma}[Bernstein Inequality]
Suppose that $(Z_i: i \geq 1)$ are mean zero, independent  rv and there exists $M>0$:
\[
|Z_i | \leq M.
\]
Then,
\[
P(\frac{1}{n} \sum_{i=1}^n Z_i >t) \leq \exp \left (- \frac{n^2 t^2/2}{\sum_{i=1}^n EZ^2_i  +M n t/3} \right ).
\]
\end{lemma}

{\bf \noindent Proof of Lemma~\ref{Lemma:1}}

Proof for (\ref{eqn_pos_part_1}) below is adapted from Bubeck et. al. 2013.

We would like to show that
\begin{equation} \label{eqn:101}
P \left (EX - \frac{1}{n}\sum_{m=1}^n X_m I(|X_m| \leq B_m) \geq
(2+\alpha)  K^{\alpha^{-1} }
\left (  \frac{\log (\delta^{-1})}{n}   \right )^{(\alpha-1)\alpha^{-1}} \right ) \leq \delta.
\end{equation}
The LHS under probability in (\ref{eqn:101}) equals
\[
\frac{1}{n} \sum_{m=1}^n ( EX - E( X I(|X| \leq B_m))   +
\frac{1}{n} \sum_{m=1}^n(E( X I(|X| \leq B_m)) - X_m I(|X_m| \leq B_m))
\]
The first term above is bounded from above by
\[
\frac{1}{n} \sum_{m=1}^n \frac{K}{B_m^{\alpha-1}}= \frac{K^{\alpha^{-1}} \log (\delta^{-1})^{(\alpha-1)\alpha^{-1}}}{n}\sum_{m=1}^n \frac{1}{m^{(\alpha-1)\alpha^{-1}}} \leq \alpha K^{\alpha^{-1}}
\left (\frac{ \log (\delta^{-1})}{n}\right )^{(\alpha-1)\alpha^{-1}}.
\]
Hence, the probability in (\ref{eqn:101}) is bounded from above by
\[
P \left (
\frac{1}{n} \sum_{m=1}^n(E( X I(|X| \leq B_m)) - X_m I(|X_m| \leq B_m))
 \geq
2  K^{\alpha^{-1} }
\left (  \frac{\log (\delta^{-1})}{n}   \right )^{(\alpha-1)\alpha^{-1}} \right )
\]
This probability can be bounded from above by $\delta$ using
Bernstein inequality above with,
$Z_m=E( X I(|X| \leq B_m)) - X_m I(|X_m| \leq B_m)$,
$M=B_n$ and
$EZ_m^2 \leq K B_n^{2-\alpha}$,
and
\[
t=2  K^{\alpha^{-1} }
\left (  \frac{\log (\delta^{-1})}{n}   \right )^{(\alpha-1)\alpha^{-1}}.
\]
Equation (\ref{eqn_pos_part_1}) follows.

To see  (\ref{eqn_pos_part_2}),  in the above analysis,  replace $X_m I(|X_m| \leq B_m)$
and $X I(|X| \leq B_m)$ by $\min(X_m,B_m)$ and $\min(X,B_m)$, respectively.
From~(\ref{eqn:error_change}), the bias term
\[
\frac{1}{n} \sum_{m=1}^n \left ( EX - E(\min(X, B_m)) \right )
\leq
\left (\frac{\alpha-1}{\alpha} \right )^{\alpha-1} \alpha^{-1}
\frac{1}{n} \sum_{m=1}^n \frac{K}{B_m^{\alpha-1}}
\leq  \left (\frac{\alpha-1}{\alpha} \right )^{\alpha-1} K^{\alpha^{-1}}
\left (\frac{ \log (\delta^{-1})}{n}\right )^{(\alpha-1)\alpha^{-1}}.
\]
The result follows by essentially the same analysis as for the truncated
rv case.
$\Box$

\vspace{0.1in}

\noindent {\bf Proof of Lemma~\ref{lemma:1123}: }  {
First note that
$t^* \leq (a+b)^2$
for if not, then
\[
(a+b) \log t^* > (a+b)^2
\]
so that
$\log t^* > (a+b)$.
Since, $ (a+b) \geq \frac{t^*}{\log t^*}$,
it follows that
$\log^2 t^* > t^*$
providing the desired contradiction (since $\log^2 t < t$ for $t \geq 1$).

Now
\[
t^* = a + b \log t^* = a + b \log (a +b \log t^*) \leq
a+ b \log (a+ 2b \log (a+b)).
\]
and the result easily follows.}
$\Box$.

\section*{REFERENCES}
\begin{hangref}

\item
Audibert, J. Y., \& Bubeck, S. 2010. Best arm identification in multi-armed bandits. In COLT-23th Conference on Learning Theory-2010 (pp. 13-p).

\item
Bechhofer, R. E., Kiefer, J., and Sobel, M. 1968. Sequential identification and ranking procedures: with special reference to Koopman-Darmois populations (Vol. 3). University of Chicago Press.
\item
Blanchet, J., Liu, J., and Zwart, B. 2008. Large deviations perspective on ordinal optimization of heavy-tailed systems. In Proceedings of the 2008 Winter Simulation Conference. 489-494. IEEE Press.

\item
Branke, J., Chick, S. E., and Schmidt, C. 2007. Selecting a selection procedure. {\em Management Science}, 53(12), 1916-1932.

\item
Broadie, M., Han, M., and Zeevi, A. 2007. Implications of heavy tails on simulation-based ordinal optimization. In Proceedings of the 2007 Winter Simulation Conference. 439-444). IEEE Press.

\item
Bubeck, S., and Cesa-Bianchi, N. 2012. Regret Analysis of Stochastic and Non-stochastic Multi-armed Bandit Problems. In Foundations and Trends in Machine Learning, Vol 5: No 1, 1-122.

\item
Bubeck, S., Cesa-Bianchi, N., and Lugosi, G. 2013. Bandits with heavy tail.  {\em IEEE Transactions on Information Theory}, 59(11), 7711-7717.

%\item
%Cesa-Bianchi, N., \& Lugosi, G. 2006. Prediction, learning, and games. Cambridge University Press.
\item
Cappe, O., Garivier, A., Maillard, O. A., Munos, R., and Stoltz, G. 2013. Kullback-Leibler upper confidence bounds for optimal sequential allocation. {\em The Annals of Statistics}, 41(3), 1516-1541.

\item Chen, C.H., J. Lin, E. Yucesan and S. E. Chick. 2000.
Simulation Budget Allocation for Further Enhancing the Efficiency
of Ordinal Optimization. {\em Journal of Discrete Event Dynamic
Systems: Theory and Applications}, {Vol. 10}, 251-270.

\item
Chernoff, H., 1959. Sequential design of experiments. {\em The Annals of Mathematical Statistics}, 30(3), pp.755-770.

\item Dai, L. 1996. Convergence Properties of Ordinal Comparison
in the Simulation of Discrete Event Dynamic Systems. {\em Journal
of Optimization Theory and Applications,} Vol. 91, 2, 363-388.

%\item
%Duffy, K., \& Metcalfe, A. P. 2005. The large deviations of estimating rate functions. Journal of applied probability, 42(1), 267-274.

\item
Even-Dar, E., Mannor, S., and Mansour, Y. 2002. PAC bounds for multi-armed bandit and Markov decision processes. Computational Learning Theory. pp. 255-270. Springer Berlin Heidelberg.

\item
Even-Dar, E., Mannor, S., and Mansour, Y. 2006. Action elimination and stopping conditions for the multi-armed bandit and reinforcement learning problems. {\em The Journal of Machine Learning Research}, 7, 1079-1105.

%\item
%Dembo, A., and O. Zeitouni. 1998. {\it Large Deviations Techniques
%and Applications.} Jones and Bartlett, Boston, MA.

%\item
%Dai, L., and Chen, C.H. 1997. Rate of convergence for ordinal comparison of
%dependent
%simulations in discrete event dynamic systems. {\em Journal of Optimization
%Theory and Applications} Vol. 94, 29-54.

\item
Foss, S., Korshunov, D., and Zachary, S. 2011. An introduction to heavy-tailed and subexponential distributions. New York: Springer.

\item
Frazier, P.I., 2014. A fully sequential elimination procedure for indifference-zone ranking and selection with tight bounds on probability of correct selection. {\em Operations Research}, 62(4), pp.926-942.

\item
Garivier, A. and Kaufmann, E., 2016. Optimal best arm identification with fixed confidence. In Conference on Learning Theory (pp. 998-1027).

\item
Glynn, P., and Juneja, S. 2004. A large deviations perspective on ordinal optimization. In Proceedings of the 2004 Winter Simulation Conference. 577-585. IEEE Press.

\item
Glynn, P. W., and Juneja, S. 2011. Ordinal optimization: A nonparametric framework. In Proceedings of the 2011 Winter Simulation Conference. 4062-4069. IEEE Press.

\item
Glynn, P. and Juneja, S. 2015. Ordinal optimization-empirical large deviations rate estimators, and stochastic multi-armed bandits. arXiv preprint arXiv:1507.04564.

\item
Glynn, P. W., and Whitt, W. 1992. The asymptotic validity of sequential stopping rules for stochastic simulations. {\em The Annals of Applied Probability}, 180-198.

\item
Glynn, P.W. and Zeevi, A., 2008. Bounding stationary expectations of Markov processes. In Markov processes and related topics: a Festschrift for Thomas G. Kurtz (pp. 195-214). Institute of Mathematical Statistics.

\item
Ho, Y. C., Sreenivas, R., and Vakili, P. 1992. Ordinal optimization of DEDS. {\em Discrete event dynamic systems}, 2(1), 61-88.

\item
Hunter, S. R., and Pasupathy, R. 2013. Optimal sampling laws for stochastically constrained simulation optimization on finite sets. {\em INFORMS Journal on Computing}, 25(3), 527-542.

\item
Jamieson, K., Malloy, M., Nowak, R. and Bubeck, S. 2014. lil?ucb: An optimal exploration algorithm for multi-armed bandits. In Conference on Learning Theory. pp. 423-439.

\item
Jennison, C., Johnstone, I.M. and Turnbull, B.W., 1982. Asymptotically optimal procedures for sequential adaptive selection of the best of several normal means. In Statistical decision theory and related topics III. pp. 55-86.

\item
Kim, S. H., and Nelson, B. L. 2001. A fully sequential procedure for indifference-zone selection in simulation. {\em ACM Transactions on Modeling and Computer Simulation}, 11(3), 251-273.

\item Kim, S. H., and B. L. Nelson. 2003. Selecting the Best
System: Theory and Methods. {\em Proceedings of the 2003 Winter
Simulation Conference,} Ed. S. Chick, P. J. Sanchez, D. Ferrin,
and D. J. Morrice. 101-112. Piscataway, New Jersey: Institute of
Electrical and Electronics Engineers.

%\item
%Kim, S. H., \& Nelson, B. L. 2006. On the asymptotic validity of fully sequential selection procedures for steady-state simulation. {\em Operations Research}, 54 (3), 475-488.

\item
Kaufmann, E., Cappe, O. and Garivier, A., 2016. On the complexity of best-arm identification in multi-armed bandit models. {\em The Journal of Machine Learning Research}, 17(1), pp.1-42.

\item
Lai, T. L., and Robbins, H. 1985. Asymptotically efficient adaptive allocation rules. Advances in applied mathematics, 6(1), 4-22.

\item
Nelson, B. L., and Matejcik, F. J. 1995. Using common random numbers for indifference-zone selection and multiple comparisons in simulation. {\em Management Science}, 41(12), 1935-1945.

%\item
%Nelson, B. L., Swann, J., Goldsman, D., \& Song, W. 2001. Simple procedures for selecting the best simulated system when the number of alternatives is large. Operations Research, 49(6), 950-963.

\item
Paulson, E., 1964. A sequential procedure for selecting the population with the largest mean from $ k $ normal populations. {\em The Annals of Mathematical Statistics}, 35(1), pp.174-180.

\item
Shin, D., Broadie, M. and Zeevi, A. 2018. Tractable Sampling Strategies for Ordinal Optimization. Preprint.

\item
Szechtman, R., and Yucesan, E. 2008. A new perspective on feasibility determination. In Proceedings of the 2008  Winter Simulation Conference.  273-28). IEEE Press.
\end{hangref}

\end{document}